\documentclass[a4paper,12pt]{amsart}
\usepackage{amssymb}
\usepackage[all,cmtip]{xy}
\usepackage{amsmath}
\usepackage{chemarrow}
\usepackage{graphicx}
\usepackage[colorlinks, %注释掉此项则交叉引用为彩色边框(将colorlinks和pdfborder同时注释掉)
pdfborder={0 0 1}, %注释掉此项则交叉引用为彩色边框
linkcolor=cyan,
anchorcolor=magenta%red,
citecolor=green
]{hyperref}
\usepackage{mathrsfs}
\usepackage{titletoc}
\usepackage{bm}
\usepackage{bbm}
\usepackage{bbold} %空心数字
\usepackage{cite}
\usepackage{subfiles}
\usepackage{extarrows} %长等号
\usepackage[hmarginratio=1:1]{geometry}
\usepackage{colordvi}
\usepackage{color}
\usepackage{stmaryrd}%\mapsfrom
\usepackage{multirow} % Required for multirows 表格里可以有多行
\usepackage{tikz-cd}
\usepackage{enumitem}
\setlist{nosep}
\makeatletter
\renewcommand{\@biblabel}[1]{[#1]\hfill}
\makeatother%文献数字编号左对齐

\allowdisplaybreaks

\DeclareMathOperator{\Aut}{Aut}

\DeclareMathOperator{\End}{End}
\DeclareMathOperator{\Ext}{Ext}

\DeclareMathOperator{\hdet}{hdet}

\DeclareMathOperator{\Hom}{Hom}

\DeclareMathOperator{\hlc}{H}

\DeclareMathOperator{\id}{id}

\DeclareMathOperator{\kk}{\Bbbk}

\DeclareMathOperator{\m}{\mathfrak{m}}

\DeclareMathOperator{\Oz}{Oz}

\DeclareMathOperator{\rank}{rank}

\newcommand{\To}{\longrightarrow}

\numberwithin{equation}{section}

\theoremstyle{definition}

\newtheorem{thm}{Theorem}[section]
\newtheorem{prop}[thm]{Proposition}
\newtheorem{lem}[thm]{Lemma}

\newtheorem{cor}[thm]{Corollary}
\newtheorem{defn}[thm]{Definition}

\newtheorem{ques}[thm]{Question}
\newtheorem{ex}[thm]{Example}

\usepackage{lineno}
\begin{document}
\author{Silu Liu}
\address{School of Mathematical Sciences, Fudan University, Shanghai 200433, China}
\email{liusl20@fudan.edu.cn}

\author{Quanshui Wu}
\address{School of Mathematical Sciences, Fudan University, Shanghai 200433, China}
\email{qswu@fudan.edu.cn}

\author{Ruipeng Zhu}
\address{School of Mathematics, Shanghai University of Finance and Economics, Shanghai 200433, China}
\email{zhuruipeng@sufe.edu.cn}

\thanks{Quanshui Wu has been supported by the NSFC 12471032, and Ruipeng Zhu has been supported by the NSFC 12301052.}

\begin{abstract}
    We prove that the ozone group of any PI Artin--Schelter regular algebra is abelian, which answers a question of Chan--Gaddis--Won--Zhang. 
    For any Calabi--Yau PI Artin--Schelter regular algebra, we  prove that the homological determinant of its ozone group acting on it is trivial. 
\end{abstract}

\title{The ozone groups of PI Artin--Schelter regular algebras are abelian}
\maketitle

%%%%%%%%%%
\section{Introduction}

The concept of ozone groups was introduced in \cite{CGWZ24, CGWZ25}. For a noncommutative ring $A$, the ozone group of $A$, denoted by $\Oz(A)$, is defined as the group of automorphisms of $A$ as a $Z(A)$-algebra, where $Z(A)$ refers to the center of $A$. According to \cite{CGWZ25} and \cite{Zhu23}, ozone groups of PI Artin--Schelter regular algebras are finite groups consisting of conjugations induced by regular normal elements. Therefore the study of ozone groups helps identifying normal elements of Artin--Schelter regular algebras. Moreover, ozone groups were applied to study the interplay between the PI Artin--Schelter regular algebras and their centers, which was effective for skew polynomial rings.

Skew polynomial rings are special examples of Artin--Schelter regular algebras. Therefore it is natural to ask whether some properties of PI skew polynomial rings are in fact shared by PI Artin--Schelter regular algebras. For instance, ozone groups of skew polynomial rings are not only finite but also abelian. Therefore a natural question was proposed in \cite[Question 0.3]{CGWZ25}.

\begin{ques}
    Is the ozone group of a noetherian PI Artin--Schelter regular algebra abelian?
\end{ques}

We answer this question affirmatively in this note.

\begin{thm} (See Theorem \ref{thm-ozone-abelian})
    Let $A$ be a noetherian PI Artin--Schelter regular algebra. Then $\Oz(A)$ is a finite abelian group.
\end{thm}

It was showed in \cite{CGWZ24} that ozone groups can be applied to characterize Calabi--Yau property of PI skew polynomial rings. To be precise, a skew polynomial ring is Calabi--Yau if and only if its ozone group action has the trivial homological determinant\cite[Theorem 0.10]{CGWZ24}. It is a widely concerned question how to characterize Calabi--Yau property of Artin--Schelter regular algebras. Consequently the following question arose.

\begin{ques}\cite[Question 5.2]{CGWZ25}
    What role does Nakayama automorphism play in the study of the ozone group and the center? If $\Oz(A)$ is assumed to be abelian, what is the connection between $A$ being Calabi--Yau and the action of $\Oz(A)$ on $A$ having trivial homological determinant?
\end{ques}

We also provide a partial result concerning this question.

\begin{thm} (See Theorem \ref{thm hdet-of-ozone})
    % Let $A$ be a noetherian Artin--Schelter regular algebra. Let $a\in A$ be a homogeneous regular normal element. Then $\mu_A(a)=\hdet_A(\eta_a)a$. Consequently, if $A$ is PI and Calabi--Yau, then the $\Oz(A)$-action on $A$ has the trivial homological determinant. \textcolor{red}{If  you write in this way, $\mu_A$ need to be defined. 这里差不多按原来的写？}
    Let $A$ be a noetherian PI Artin--Schelter regular algebra which is Calabi--Yau. Then the action of $\Oz(A)$ on $A$ has the trivial homological determinant.
\end{thm}

%%%%%%%%%%
\section{Preliminaries}
Throughout $\kk$ is assumed to be an algebraically closed field of characteristic $0$. For a ring $A$, $A^o$ refers to the opposite ring of $A$. By $A$-modules we mean left $A$-modules.

A \textbf{polynomial identity ring} (or PI ring for short) is a ring satisfying a polynomial identity. A ring which is finitely generated as a module over a central subring is a PI ring. The converse  does not hold in general. We refer to \cite[Chapter 13]{MR01} for the facts about PI rings. The main fact we need about PI rings is the following structure theorem.

\begin{prop}(Posner's theorem, \cite[Theorem 6.5]{MR01})\label{prop Posner's thm}
    Let $A$ be a prime PI ring with center $Z$. Let $K$ be the fraction field of $Z$, and let $D:=A\otimes_Z K$. Then $D$ is a central simple $K$-algebra.
\end{prop}

A $\kk$-algebra $A$ is \textbf{connected graded} if $A$ is an $\mathbb{N}$-graded algebra with $A_0=\kk$.

\begin{defn}
    A noetherian connected graded $\kk$-algebra $A$ is called \textbf{Artin--Schelter Gorenstein} (or AS Gorenstein for short) of dimension $d$ if 
    \begin{enumerate}
        \item $A$ has finite left and right injective dimension;
        \item 
        $\begin{cases}
          \Ext_A^i(\kk,A)\cong\Ext_{A^{o}}^i(\kk,A)=0, \forall\ i\neq d;\\
        \Ext_A^d(\kk,A)\cong\Ext_{A^{o}}^d(\kk,A)\cong \kk(\mathfrak{l}) \text{ for some } \mathfrak{l}\in \mathbb{Z},
        \end{cases}$ 
    \end{enumerate} 
     where $\mathfrak{l}$ is called the \textbf{Gorenstein parameter} of $A$.
     
    If moreover $A$ has finite global dimension then $A$ is called \textbf{Artin--Schelter regular} (or AS regular for short). 
\end{defn}

Noetherian PI AS regular algebras were well studied in \cite{SZ94}.

\begin{thm}\cite[Corollary 1.2]{SZ94}\label{thm AS-domain}
    Let $A$ be a noetherian PI AS regular algebra. Then $A$ is a domain and a finitely generated module over its center $Z$, which is a Krull domain.
\end{thm}

Homological determinant was introduced by J\o rgensen and Zhang in \cite{JZ00}. Assume that $A$ is AS Gorenstein of dimension $d$ with Gorenstein parameter $\mathfrak{l}$. Let $\sigma\in\Aut_{gr}(A)$, the group of graded $\kk$-algebra automorphisms of $A$. For $A$-modules $M$ and $N$, a $\kk$-linear morphism $f:M\to N$ is called \textbf{$\sigma$-linear} if $f(am)=\sigma(a)f(m)$ for all $m\in M$ and $a\in A$. 

Let $(-)'$ be the graded duality functor $\underline{\Hom}_{\kk}(-,\kk)$. By \cite[Lemma 2.2]{JZ00} $\sigma: A\to A$, $a\mapsto \sigma(a)$ is a $\sigma$-linear morphism between graded $A$-modules, which induces a $\sigma$-linear morphism $\hlc_{\m_A}^d(\sigma): \hlc_{\m_A}^d(A)\to \hlc_{\m_A}^d(A)$ between the local cohomologies. Since $A$ is AS Gorenstein, $\hlc_{\m_A}^d(A)\cong A'(\mathfrak{l})$\cite[Lemma 2.1]{JZ00}. Then there exists unique $c_{\sigma}\in\kk^{\times}$ such that $\hlc_{\m_A}^d(\sigma)=c_{\sigma}(\sigma^{-1})'$ \cite[Lemma 2.2]{JZ00}. 

\begin{defn}\cite[Definition 2.3]{JZ00}
    The constant $c_{\sigma}^{-1}$ above is called the \textbf{homological determinant} of $\sigma$ acting on $A$, which is denoted as $\hdet_A(\sigma)=c_{\sigma}^{-1}$. 
\end{defn}

Consider the map $\hdet_A: \Aut_{gr}(A)  \to \kk^{\times}$, $\sigma\mapsto c_{\sigma}^{-1}$. By \cite[Proposition 2.5]{JZ00} $\hdet_A(\sigma\tau)=\hdet_A(\sigma)\hdet_A(\tau)$ for all $\sigma,\tau\in \Aut_{gr}(A)$. If $\sigma$ is of finite order, then it follows that $\hdet_A(\sigma)$ is a root of unity. For any $G \leqslant \Aut_{gr}(A)$, if $\hdet_A(\sigma)=1$ for all $\sigma \in G$, then we say that the $G$-action on $A$ has the \textbf{trivial homological determinant}.

%%%%%%%%%%
\section{Ozone groups of PI AS regular algebras}

\begin{defn}\cite[Definition 1.1]{CGWZ25}
    Let $A$ be a ring with center $Z$. The \textbf{ozone group} of $A$ is defined as
    $\Oz(A):=\{\sigma\in\Aut(A)\mid \sigma(z)=z, \forall z\in Z\}$.
\end{defn}

An element $a$ of a ring $A$ is called \textbf{normal} if $Aa=aA$. If $ab\neq 0$ and $ba\neq 0$ for any nonzero element $b$ of $A$, then $a$ is called a \textbf{regular} element. Note that normal elements in a prime ring are always regular. 

Let $a$ be a regular normal element of $A$. Then $a$ induces an automorphism $\eta_a$ of $A$, given by $ax=\eta_a(x)a$, for all $x\in A$. It is straightforward to verify that, for regular normal elements $a,b\in A$, the product $ab$ is a regular normal element of $A$ and $\eta_{ab}=\eta_a\eta_b$. Obviously, $\eta_a\in\Oz(A)$. In fact, under some mild assumptions, the elements in $\Oz(A)$ are of the above form. 

\begin{lem} \cite[Lemmas 1.9 and 1.10]{CGWZ25}\label{lem ozone-normal-elements}
    Let $A$ be a prime ring which is a finitely generated module over its center $Z$. Let $\sigma\in\Oz(A)$.
    \begin{enumerate}
        \item There exists a nonzero normal element $a\in A$ such that $ab=\sigma(b)a$ for all $b\in A$. In other words, $\sigma=\eta_a$.
        \item If further $A$ is a $\mathbb{Z}$-graded domain then $\sigma=\eta_a$ for some nonzero homogeneous normal element $a\in A$. Consequently $\Oz(A)\leqslant \Aut_{gr}(A)$.
        \item For any nonzero normal elements $a,b\in A$, $\eta_a=\eta_b$ if and only if there exist nonzero central elements $z,z'$ of $A$ such that $za=z'b$.
    \end{enumerate}
\end{lem}

\begin{lem}\label{lem ozone-abelian}
    Let $A$ be a domain such that $A$ is a finitely generated module over its center $Z$, which contains all the roots of unity.
    Suppose $a$ and $b$ are two nonzero normal elements of $A$. If the order of the automorphism $\eta_a\eta_b\eta_a^{-1}\eta_b^{-1}$ is finite, then $ab = \zeta ba$ for some root of unity $\zeta$.
    % Let $A$ be a $\kk$-algebra which is a domain and a finitely generated module over its center $Z$. Let $a,b$ be nonzero normal elements of $A$.
    % If the automorphism $\eta_a\eta_b\eta_a^{-1}\eta_b^{-1}$ of $A$ is of finite order, then $ab = \zeta ba$ for some root of unity $\zeta$.
\end{lem}

\begin{proof}
    Let $K$ be the fraction field of $Z$. Then, by Proposition \ref{prop Posner's thm}, $D:=A\otimes_ZK$ is a central simple $K$-algebra. Moreover, since $A$ is a domain, $D$ is a domain. By the Artin--Wedderburn theorem, $D$ is a division ring. So, $a$ and $b$ are invertible elements in $D$.
    
    Note that for any $\sigma\in\Oz(A)$, $\sigma$ can be extended uniquely to an automorphism of $D$ by $\sigma(xz^{-1}):=\sigma(x)z^{-1}$ for all $x\in A$ and $z\in Z$. Let $u:=aba^{-1}b^{-1} \in D$. Then
    \[ ux = (aba^{-1}b^{-1})x = \eta_{aba^{-1}b^{-1}}(x)(aba^{-1}b^{-1})= \big(\eta_a\eta_b\eta_a^{-1}\eta_b^{-1}(x)\big)u, \ \forall\ x \in D.\]
    If $(\eta_a\eta_b\eta_a^{-1}\eta_b^{-1})^m=\id_A$, then $u^mx= xu^m$ for all $x \in D$. Hence $u^m \in Z(D) = K$.

    Consider the following norm map $N:D \To K$, which is the composition
    \[D \xrightarrow{\ell} \End_K(D)=M_n(K)\xrightarrow{\det}K\]
    where $\ell: D \to \End_K(D)$ is given by the left multiplications.
    Since the norm map $N:D \To K$ is multiplicative, it follows that $$N(u)=N(a)N(b)N(a^{-1})N(b^{-1})=1, \textrm{ and } 1=N(u)^m=N(u^m)=(u^m)^n$$ as $u^m\in K$.
    
    By assumption, $Z$ contains all $mn$-th roots of unity ${\zeta_1,\cdots, \zeta_{mn}}$.
    Then 
    $$t^{mn} - 1 = \prod_{i=1}^{mn} (t - \zeta_i) \in Z[t] \subseteq K[t].$$ 
    It follows from $u^{mn} = 1$ that $\prod_{i=1}^{mn} (u - \zeta_i) = 0 \in D$. Since $D$ is a domain, $u = \zeta_i \in Z \subseteq K$ for some $i$.
    Hence $ab = \zeta_i ba$.
\end{proof}

Note that if $A$ is a $\kk$-algebra, then its center automatically contains all roots of unity, since $\kk$ is assumed to be algebraically closed. By Theorem \ref{thm AS-domain}, Lemma \ref{lem ozone-normal-elements} and Lemma \ref{lem ozone-abelian} are applicable to noetherian PI AS regular algebras.

\begin{thm}\label{thm-ozone-abelian}
    Let $A$ be a noetherian PI AS regular algebra. Then $\Oz(A)$ is a finite abelian group.
\end{thm}

\begin{proof} 
    The finiteness of $\Oz(A)$ is proved in \cite[Theorem 2.5]{Zhu23} and \cite[Theorem E]{CGWZ25}.
    It suffices to show that $\sigma\tau = \tau\sigma$ for any $\sigma, \tau \in \Oz(A)$.
    By Lemma \ref{lem ozone-normal-elements} (1), there exist nonzero normal elements $a,b\in A$ such that $\sigma=\eta_a$, $\tau=\eta_b$.
    Lemma \ref{lem ozone-abelian} then guarantees the existence of a root of unity $\zeta$ such that $ab = \zeta ba$.
    By Lemma \ref{lem ozone-normal-elements} (3)
    \[\sigma\tau=\eta_a\eta_b=\eta_{ab}=\eta_{\zeta ba}=\eta_{ba}=\eta_b\eta_a=\tau\sigma,\]
    which completes the proof.
    % Since $x\in A_{\sigma}$ and $y \in A_{\tau}$, the product $xy$ lies in $A_{\sigma \tau}$. Similarly, $yx \in A_{\tau \sigma}$.
    % Then $0 \neq xy = \zeta yx$ implies $A_{\sigma\tau} \cap A_{\tau\sigma}\neq 0$, which forces $\sigma\tau = \tau\sigma$ by Lemma \ref{lem ozone-normal-elements}.
\end{proof}

As a result of Theorem \ref{thm-ozone-abelian}, some results of \cite{CGWZ25} can be restated as follows.

\begin{cor}\label{cor restatement}
    Let $A$ be a noetherian PI AS regular algebra with center $Z$.
    \begin{enumerate}
        \item \cite[Theorem B]{CGWZ25} A finite group $G$ is isomorphic to the ozone group of some noetherian PI AS regular algebra if and only if $G$ is abelian.
        \item \cite[Lemma 4.10]{CGWZ25} $\sum_{\sigma\in\Oz(A)}A_{\sigma}$ is a direct sum of $Z$-modules, where $A_{\sigma}:=\{a\in A\mid \sigma(x)a=ax, \forall x\in A\}$.
        \item \cite[Theorem F]{CGWZ25} Suppose moreover that $A$ is generated in degree one. Then the following are equivalent.
        \begin{enumerate}
            \item $A$ is a skew polynomial ring.
            \item $|\Oz(A)|=\rank_ZA$.
            \item $A$ is generated by normal elements.
            \item $A^{\Oz(A)}=Z$.
        \end{enumerate}
    \end{enumerate}
\end{cor}

Although the ozone group action of $A$ may be complicated, the action of $\Oz(A)$ on the normal elements of $A$ is quite simple.

\begin{cor}\label{cor rou}
    Let $a$ be a normal element of a noetherian PI AS regular algebra $A$. For any $\tau \in \Oz(A)$, $\tau(a)=\zeta a$ for some root of unity $\zeta$.
\end{cor}

\begin{proof}
    We may assume that $a\neq 0$. Take $0\neq b\in A_{\tau}$. It follows from Lemma \ref{lem ozone-abelian} that there exists some root of unity $\zeta$ such that $ab=\zeta ba$. By the definition of $A_{\tau}$, $ba= \tau(a)b$. The conclusion follows as $A$ is a domain.
\end{proof}

An AS regular algebra $A$ is called \textbf{Calabi--Yau} if its Nakayama automorphism $\mu_A$ of $A$ is the identity map. Note that $\mu_A\in\Oz(A)$\cite[Proposition 4.4]{BZ08}. The homological determinants of the $\Oz(A)$-action on $A$ can be calculated through $\mu_A$.

\begin{lem}\cite[Corollary 5.4]{RRZ17}\label{lemma hdet-of-nak}
    Let $A$ be an AS Gorenstein algebra which is a homomorphic image of a noetherian AS regular algebra. Then $\hdet_A(\mu_A)=1$.
\end{lem}

Recall that for any $\sigma\in\Aut(A)$, an element $a\in A$ is called \textbf{$\sigma$-skew} if $ax=\sigma(x)a$ for all $x\in A$ and $\sigma(a)=a$ \cite[$\S$ 2]{ZSL20}. 

\begin{thm}\label{thm hdet-of-ozone}
    Let $A$ be a noetherian AS regular algebra. Let $a\in A$ be a homogeneous regular normal element. Then $\mu_A(a)=\hdet_A(\eta_a)a$. Consequently, if $A$ is PI and Calabi--Yau, then the $\Oz(A)$-action on $A$ has the trivial homological determinant.
\end{thm}

\begin{proof}
    We may assume that $a$ is of positive degree. By definition, $a$ is an $\eta_a$-skew element. It follows from \cite[Proposition 2.3]{ZSL20} that $\mu_A(a)=\zeta a$ for some $\zeta\in\kk^{\times}$. By Rees' lemma $\Bar{A}:=A/(a)$ is AS Gorenstein. It follows from \cite[Lemma 1.5]{RRZ14} that 
    \begin{align*}
        \mu_{\Bar{A}}=(\mu_A\circ\eta_a)|_{\Bar{A}}.
    \end{align*}
    Note that $\eta_a(a)=a$. According to \cite[Proposition 2.4]{JZ99},
    \begin{align*}
        \hdet_A(\eta_a)=&\hdet_{\Bar{A}}(\eta_a|_{\Bar{A}});\\
        \hdet_A(\mu_A)=&\zeta\hdet_{\Bar{A}}(\mu_A|_{\Bar{A}}).
    \end{align*}
    By Lemma \ref{lemma hdet-of-nak},
    \[ 1=\hdet_{\Bar{A}}(\mu_{\Bar{A}})=\hdet_{\Bar{A}}(\mu_A|_{\Bar{A}})\cdot \hdet_{\Bar{A}}(\eta_a|_{\Bar{A}})=\zeta^{-1}\cdot\hdet_A(\eta_a),\]
    and the rest follows immediately.
\end{proof}

% \begin{thm}\label{thm-hdet-of-ozone}
%     Let $A$ be a noetherian PI AS regular algebra. For any $\sigma\in\Oz(A)$ and $a\in A_{\sigma}$, $\mu_A(a)=\hdet_A(\sigma)a$. Consequently, if $A$ is moreover Calabi--Yau, then the $\Oz(A)$-action on $A$ has the trivial homological determinant.
% \end{thm}

% \begin{proof}
%    Obviously, $\sigma(a)=a$ as $A$ is a domain. By Corollary \ref{cor rou}, $\mu_A(a)=\zeta a$ for some root of unity $\zeta$. We may assume $a\neq 0$. By Rees' lemma $A/(a)$ is AS Gorenstein. It follows from \cite[Lemma 1.5]{RRZ14} that 
%    \begin{align*}
%        \mu_{A/(a)}=(\mu_A\circ\sigma)|_{A/(a)}.
%    \end{align*}
%     According to \cite[Proposition 2.4]{JZ99},
%     \begin{align*}
%         \hdet_A(\sigma)=&\hdet_{A/(a)}(\sigma|_{A/(a)});\\
%         \hdet_A(\mu_A)=&\zeta\hdet_{A/(a)}(\mu_A|_{A/(a)}).
%     \end{align*}
%    The rest follows from Lemma \ref{lemma-hdet-of-nak}.
% \end{proof}

From Theorem \ref{thm hdet-of-ozone}, we conclude that for any noetherian PI AS regular algebra $A$, $\Oz(A)$ acts on $A$ with the trivial homological determinant if and only if $\mu_A$ fixes all normal elements. Since skew polynomial rings are generated by normal elements, $\Oz(A)$ acting on $A$ with the trivial homological determinant indicates that $\mu_A$ fixes $A$, or equivalently, that $A$ is Calabi--Yau. See \cite[Theorem 0.10]{CGWZ24} for further connections. At the end, we show by an example that this  does not necessarily hold for PI AS regular algebras.

\begin{ex}\label{ex trivial-det-not-CY}
    Consider the down-up algebra
    $$A:=A(0,1)=\kk\langle x,y\rangle/(x^2y-yx^2, xy^2-y^2x).$$
    Then $\Oz(A)=\{1,\mu_A\}$, where $\mu_A:x\mapsto -x, y\mapsto -y$\cite[Example 5.4]{CGWZ25}. Therefore, $A$ is not Calabi--Yau. By Lemma \ref{lemma hdet-of-nak} $\hdet_A(\mu_A)=1$, which implies that $\hdet_A(\Oz(A))=1$.

\end{ex}

\section*{Acknowledgments}
The authors express their deep appreciation for the invaluable comments and suggestions provided by the referee.

\nolinenumbers

\bibliographystyle{amsalpha}
\bibliography{refs.bib}

\end{document}